\documentclass[11pt]{article}

\usepackage{amsmath,latexsym,epsfig}
\usepackage{amsfonts,amssymb,amscd}

\title{\bf The  number of ramified coverings of the sphere by the  double
torus, and a general form for higher genera}

\author{
I.P.Goulden\thanks{Dept. of Combinatorics and Optimization,
University of Waterloo, Waterloo, Ontario, Canada} and
D.M.Jackson\thanks{Dept. of Combinatorics and Optimization,
University of Waterloo, Waterloo, Ontario, Canada} 
}
\date{January 15, 1999}

\begin{document}
 \maketitle

 \newtheorem{theorem}{Theorem}[section]
 \newtheorem{proposition}[theorem]{Proposition}
 \newtheorem{definition}[theorem]{Definition}
 \newtheorem{axiom}[theorem]{Axiom}
 \newtheorem{lemma}[theorem]{Lemma}
 \newtheorem{corollary}[theorem]{Corollary}
 \newtheorem{remark}[theorem]{Remark}
 \newtheorem{example}[theorem]{Example}
 \newtheorem{conjecture}[theorem]{Conjecture}

\def\sP{{\sf{P}}}
\def\sX{{\sf{X}}}
\def\bfp{{\rm\bf p}}
\def\cC{{\cal{C}}}
\def\cH{{\cal{H}}}
\def\cP{{\cal{P}}}
\def\cU{{\cal{U}}}
\def\symgp{{\mathfrak{S}}}
\def\rats{{\mathbb{Q}}}
\def\proof{{\rm\bf Proof:\quad}}
\def\qed{{\hfill{\large $\Box$}}}
\def\sph{\mathbb{S}^2}
\def\surf{\sX}
\def\pdif#1#2{\frac{\partial{#1}}{\partial p_{#2}}}
\def\dpdif#1#2#3{\frac{\partial^2{#1}}{\partial p_{#2}\partial p_{#3}}}
\def\macz#1{\vartheta(#1)}
\def\wij#1#2{w^{(#1)}_{#2}}

\def\gp#1#2{\frac{c_#1(1^#2)}{#2!}}
\def\dx{D}
\def\dxop#1{x\frac{\partial #1}{\partial x}}
\def\ddx#1{D^#1}

\def\atp#1#2{\stackrel{\scriptstyle{#1}}{\scriptstyle{#2}}}

\begin{abstract}
An explicit expression is obtained for the generating series for the
number of ramified coverings of the sphere by the double torus, with elementary
branch points and prescribed ramification type over infinity.  
Thus we are able to determine various linear recurrence equations
for the  numbers of these coverings with no ramification over infinity;
one of these recurrence equations has previously been conjectured by
Graber and Pandharipande.
The general form of this series is conjectured for the number of these coverings by a surface of
arbitrary genus that is at least two.
\end{abstract}


\section{Introduction and Background}\label{SIn}
Let $\mu_m^{(g)}(\alpha)$ be the number of almost simple ramified coverings
of $\sph$ by $\surf$ with ramification type $\alpha$ where $\surf$ is a 
compact connected Riemann surface of genus $g,$ and $\alpha$ is a partition with $m$ parts.
The problem of determining an (explicit) expression for $\mu_m^{(g)}(\alpha)$ is called the
{\em Hurwitz Enumeration Problem},  a brief account of which is given in~\cite{GJdec11}.
The terminology and notation used here will be consistent with the latter paper.

There appears to be a natural topological distinction between the low genera instances
of the problem, namely for
$g\le1,$ on the one hand, and the higher genera instances, namely for $g\ge2,$ on the
other hand. In this paper
we address the higher genera case of the Hurwitz Enumeration Problem.
The distinction between the low genera and the high genera cases manifests itself in this
paper in the fact
that a general form can be given for the higher genera case, but that does not specialize to the
low genera case.
In this paper we prove an explicit result for g=2, the double torus.
Thus we are able to determine various linear recurrence equations for
$\mu_m^{(2)}(1^m)$,
corresponding to the case of no ramification over infinity.
One of these recurrence equations has previously been
conjectured by Graber and Pandharipande~\cite{V1}.
Moreover, we conjecture an explicit result for $g=3.$ 
We also give a conjecture for the general form of the generating series for $\mu_m^{(g)}(\alpha)$
for arbitrary $g\ge2.$

Let $\cC_\alpha$ be the conjugacy class of the symmetric group
$\symgp_n$ on $n$ symbols indexed by the partition $\alpha$ of $n.$ 
Let $\macz\alpha=\prod_{i\ge1} i^{m_i} m_i!,$  where $\alpha$ has $m_i$ parts equal to $i$ where $i\ge1.$
We write $\alpha\vdash n$ to indicate that
$\alpha=(\alpha_1,\ldots,\alpha_r)$ is a partition of $n$ and 
$\alpha\models n$ to indicate that $\alpha$ is a partition
of $n$ with no part equal to  one. The length, $r,$ of $\alpha$ is denoted by $l(\alpha).$

Let the generating series for $\mu_m^{(g)}(\alpha)$ be defined by
\begin{eqnarray}\label{e1}
F_g(x,\bfp) = \sum_{m,n\ge1}
\sum_{\atp{\alpha\vdash n} {l(\alpha)=m}} \frac{\mu_m^{(g)}(\alpha)}{(n+m+2g-2)!}
{p_\alpha} x^n.
\end{eqnarray}
where $p_i$, for $i\ge1,$ are indeterminates and $p_\alpha=p_{\alpha_1}p_{\alpha_2}\ldots.$
Let
\begin{eqnarray}\label{e2}
\psi_i(x,\bfp)=\sum_{n\ge1}n^{i-1}a_np_nx^n
\end{eqnarray}
where $a_n=n^n/(n-1)!$ for $n\ge1.$
Then, as we have shown in previous work~(\cite{GJtransf},\cite{GJdec11})
 $F_0$ and $F_1$ can be expressed succinctly in terms of $\psi_i\equiv\psi_i(s,\bfp)$
where $s\equiv s(x,\bfp)$ is the unique solution of the functional equation
\begin{eqnarray}\label{e3}
s=x e^{\psi_0(s,\bfp)}.
\end{eqnarray}
From~\cite{GJtransf} (Prop.~3.1) and~\cite{GJdec11} (Thm~{4.2}),
the expressions for $F_0$ and $F_1$ in terms of the $\psi_i$ are given by
\begin{eqnarray}
\left( x\frac{\partial}{\partial x} \right)^2 F_0(x,\bfp) &=& \psi_0, \label{e4}\\
F_1(x,\bfp) &=& \frac{1}{24}\left(\log(1-\psi_1)^{-1}-\psi_0\right). \label{e4a}
\end{eqnarray}

In  this paper we prove the following explicit expression for $F_2.$
\begin{theorem}\label{T1}
\begin{eqnarray}\label{3a}
F_2(x,\bfp)= \frac1{5760}
\left(
\frac{Q_3}{(1-\psi_1)^3} + \frac{Q_4}{(1-\psi_1)^4} + \frac{Q_5}{(1-\psi_1)^5}\right),
\end{eqnarray}
where
$Q_3=5\psi_4-12\psi_3+7\psi_2,$\quad $Q_4=29\psi_2\psi_3-25\psi_2^2,$ \quad and \quad  $Q_5=28\psi_2^3.$
\end{theorem}
The form of this theorem was arrived at through a careful examination and further analysis of the expressions
for $\mu_m^{(2)}(\alpha)$ for $m\le3$ that appear in the Appendix of~\cite{GJV}.
In addition, we make the following conjecture for the general form of $F_g.$
\begin{conjecture}\label{Con1}
For $g\ge2,$
$$
F_g(x,\bfp)=\sum_{d=2g-1}^{5g-5} \frac{1}{(1-\psi_1)^d}
\sum_{n=d-1}^{d+g-1}
\sum_{\atp{\theta\models n} {l(\theta)=d-2(g-1)}}
K_\theta^{(g)}\psi_\theta.
$$
where $K_\theta^{(g)}\in\rats,$ and $\psi_\theta=\psi_{\theta_1}\psi_{\theta_2}\ldots.$
\end{conjecture}

\medskip

This expression for $F_g(x,\bfp)$ is a sum of rational functions of the $\psi_i$'s with particularly simple denominators and
with numerators of prescribed form.
For $g=3$ we have determined the $K_\alpha^{(g)}$explicitly, based  on this form, with the aid of {\sf Maple}.
The resulting expression for $F_3(x,\bfp)$ is displayed in Appendix~\ref{A0}.

For $l(\theta)=1$ and for all $g,$ $K_\theta^{(g)}$ may be obtained quite readily as follows.
From~\cite{JA},
$$\frac{n\mu_1^{(g)}(n)}{(n+2g-1)!} = a_n\frac{n^{2g-2}}{2^{2g}}[x^{2g}]\left(\frac{\sinh(x)}{x}\right)^{n-1},$$
where $[x^{2g}]$ denotes the operator giving the coefficient of $x^{2g}$ in a formal power series.
But, as a polynomial in $n,$
$$[x^{2g}]\left(\frac{\sinh(x)}{x}\right)^{n-1}=
\frac{1}{6^gg!}n^g+\cdots+ \frac{2(1-2^{2g-1})B_{2g} }{(2g)!}n^0$$
where $B_{2g}$ is a Bernoulli number. 
Thus, from the conjectured form, we have
$$F_g =
\frac{\frac{1}{24^gg!}\psi_{3g-2}+\cdots+ \frac{(1-2^{2g-1})B_{2g}}{2^{2g-1}(2g)!}\psi_{2g-2}}{(1-\psi)^{2g-1}}+\cdots+
 \frac{A_g\psi_2^{3g-3}}{(1-\psi)^{5g-5}}.
$$
where $A_g$ is  a rational number.

\section{Proof of the main result}\label{Stde}
We use the approach that has been developed in~\cite{GJtransf} and extended in~\cite{GJdec11} that
makes use of  Hurwitz's encoding~\cite{Hrfgv} of the problem as a transitive ordered factorization of
a permutation with prescribed cycle type into transpositions. It has already been shown~\cite{GJV}  
by a cut-and-join analysis of the action of a transposition on the cycle structure of an arbitrary
permutation in $\symgp_n,$
that if
$$F=\widetilde{F}_0+z\widetilde{F}_1+z^2\widetilde{F}_2+\cdots$$
where $\widetilde{F}_g$ is obtained from $F_g$ by multiplying the summand by $u^{n+m+2g-2},$
then $F$
satisfies the partial differential equation
\begin{eqnarray}\label{e21}
\frac{\partial F}{\partial u} = \frac12 \sum_{i,j\ge1} \left(
ijp_{i+j}z\dpdif Fij + ijp_{i+j}
\pdif Fi \pdif Fj +(i+j)p_ip_j \pdif F{i+j} \right).
\end{eqnarray}
The techniques developed in~\cite{GJtransf} and in~\cite{GJdec11} enable us to confirm whether
a series satisfies the partial differential equation induced from this by considering only terms
of given degree in $z,$ thus grading by genus, but we do not yet possess
a method for constructing the solution of such an equation in closed form.
The next result gives the linear first order partial differential equation for $F_2$
that is induced by restricting~(\ref{e21}) in this way to terms of degree exactly two in $z$.

\begin{lemma}\label{L1}
The series $f=F_2$ satisfies the partial differential equation
\begin{eqnarray}\label{e21a}
T_0 f- T_1=0
\end{eqnarray}
where
\begin{eqnarray*}
T_0 &=& x\frac{\partial}{\partial x} + \sum_{i\ge1} p_i\pdif {}i +2
- \sum_{i,j\ge1} ij p_{i+j} \pdif{F_0}i \pdif{}j
-\frac12 \sum_{i,j\ge1} (i+j) p_ip_j \pdif{}{i+j}, \\
T_1 &=& \frac12 \sum_{i,j\ge1}  ij p_{i+j} \pdif{F_1}i \pdif{F_1}j + \frac12 \sum_{i,j\ge1} ij p_{i+j} \dpdif{F_1}ij.
\end{eqnarray*}
\end{lemma}
\proof Clearly, 
$$u\frac{\partial}{\partial u} [z^2]\, F = 
\left( 
x\frac{\partial}{\partial x} + \sum_{i\ge1} p_i\pdif {}i +2
\right) 
[z^2]\,F,$$
and the result follows by applying $[z^2]$ to~(\ref{e21}). \qed

\medskip

Lemma~\ref{L1} gives a linear partial differential equation for $F_2$ with coefficients that involve
the known series $F_0$ and $F_1$ (see~(\ref{e4}) and~(\ref{e4a})).
The proof of Theorem~\ref{T1} consists of showing that the expression for $F_2$ given in~(\ref{3a}) does indeed
satisfy this linear partial differential equation. 
To establish this, extensive use will be made of a number of results in~\cite{GJdec11} that were
used in the determination of the generating series  $F_1,$ the case of the torus. These 
enable us to reduce $T_0F_2-T_1$ to a polynomial in a new set of variables,
and this is shown to be identically zero, with the aid of {\sf Maple} to carry out the substantial
amount of routine simplification. The details of this part of the proof are suppressed, but
enough information is retained to permit it to be reproduced.

In particular
we require the following results from~\cite{GJdec11}, which are included for completeness.
They follow from~(\ref{e2}), (\ref{e3}), (\ref{e4}) and (\ref{e4a}), with the aid of Lagrange's
Implicit Function Theorem.
\begin{eqnarray}
x\frac{\partial s}{\partial x} =  \frac{s}{1-\psi_1},  \quad
\pdif sk =\frac1k \frac{a_k s^{k+1}}{1-\psi_1}, \label{e22} \label{e22a}
\end{eqnarray}
\begin{eqnarray}
x\frac{\partial \psi_i}{\partial x} =  \frac{\psi_{i+1}}{1-\psi_1},  \quad 
\pdif {\psi_i}k=  k^{i-1}a_ks^k + \frac{a_k}{k} \frac{\psi_{i+1}s^k}{1-\psi_1},\label{e23}\label{e23a}  
\end{eqnarray}
and
\begin{eqnarray}
\pdif {F_0}k &=& \frac{a_k}{k^3}s^k -\frac{a_k}{k^2}\sum_{r\ge1}  a_r p_r\frac{s^{k+r}}{k+r}, \label{e25} \\
\pdif{F_1}k &=&
\frac1{24} \frac{a_k}{k}s^k \left( \frac{1-k}{1-\psi_1} + \frac{\psi_2}{(1-\psi_1)^2}\right).   \label{e26}
\end{eqnarray}
 
It follows immediately from~(\ref{e26}) that
\begin{eqnarray}
\dpdif{F_1}ij =\frac{1}{24} \frac{a_is^i}{i}\,\frac{a_js^j}{j}
\left(
\frac{i^2+j^2+ij-i-j}{(1-\psi_1)^2} + \frac{2(i+j)\psi_2-\psi_2+\psi_3}{(1-\psi_1)^3}
+\frac{2\psi_2^2}{(1-\psi_1)^4}\right). \label{e27}
\end{eqnarray}
These account for all of the terms in the partial differential equation~(\ref{e21a})
that do not involve $F_2.$  

\medskip

\noindent\proof {\bf [Theorem~\ref{T1}]}
Let $G_2$ denote the series written explicitly on the right hand side of~(\ref{3a}).
To prove Theorem~\ref{T1} it is sufficient to  show that $G_2$ is a solution of~(\ref{e21a}),
since $G_2,$ with constant term equal to zero, clearly satisfies the initial condition for $F_2$. 
The requisite partial derivatives of $G_2$ are obtained
indirectly from~(\ref{e23}) by differentiating with respect to the $\psi_i$'s,
giving
\begin{eqnarray}
x\frac{\partial G_2}{\partial x} &=& \sum_{m\ge1}\left(\frac{\partial G_2}{\partial \psi_m}\right)
\left(x\frac{\partial\psi_m}{\partial x}\right)
=\sum_{m\ge1} \frac{\psi_{m+1}}{1-\psi_1}\frac{\partial G_2}{\partial\psi_m}, \label{e28} \\
\pdif{G_2}k &=& \sum_{m\ge1} \frac{a_ks^k}{k}\left(k^m +\frac{\psi_{m+1}}{1-\psi_1}\right)
\frac{\partial G_2}{\partial \psi_m}. \label{e28a}
\end{eqnarray}

By inspection, each summand of $T_0G_2-T_1$ may be expressed through~(\ref{e25}), (\ref{e26}), (\ref{e27}),
(\ref{e28}) and~(\ref{e28a}) 
entirely in terms of $q_i$'s where
$q_i=s^ip_i$ for $i\ge1.$ The $q_i$'s are algebraically independent, and the dependency on $s$
is perfectly subsumed. Now, in terms of these $q_i$'s, let
\begin{eqnarray}\label{e29}
M_{k,l} = \sum_{i,j\ge1} q_{i+j} a_i a_j  i^k j^l,
\end{eqnarray}
so $M_{k,l}=M_{l,k},$ and $M_{k,l}$ is homogeneous of degree one in the $q_i$'s.
Let
\begin{eqnarray} 
N_k &=& \frac12 \sum_{i,j\ge1} \left( q_iq_j a_{i+j}(i+j)^k
-2q_{i+j}  a_j j^k \frac{a_i}i \sum_{r\ge1} q_r\frac{a_r}{i+r}\right), \label{e211}
\end{eqnarray}
so $N_k$ is homogeneous of degree two in the $q_i$'s.
Then, from~(\ref{e26}) and~(\ref{e27})
\begin{eqnarray}\label{e212}
T_1 = \frac{1}{5760}\sum_{i=2}^4 \frac{S_i}{(1-\psi_1)^i},
\end{eqnarray}
where
\begin{eqnarray*}
S_2 &=& 240 M_{2,0}+125 M_{1,1} - 250 M_{1,0} + 5 M_{0,0}, \\
S_3 &=& \psi_2\left( 490 M_{1,0} -130 M_{0,0} \right) + 120 \psi_3 M_{0,0}, \\
S_4 &=& 245 M_{0,0},
\end{eqnarray*}
and, from~(\ref{e25}), (\ref{e28}), (\ref{e28a}), (\ref{e29}) and (\ref{e211}), 
\begin{eqnarray} \label{e213}
T_0 G_2 = 2G_2+\sum_{m\ge1} 
\left(
\frac{\psi_{m+1}}{1-\psi_1} \left( 1+\psi_0-M_{-2,0}-N_0 \right) +\psi_m-M_{-2,m}-N_m
\right) \frac{\partial G_2}{\partial\psi_m}.
\end{eqnarray}
But
$$
 \frac{\partial G_2}{\partial\psi_m} = \frac{1}{5760}\sum_{l} \frac{R_{l,m}}{(1-\psi_1)^l},
$$
where
$R_{l,1}=(l-1)Q_{l-1}$ for $l=4,5,6,$ $R_{3,2}=7,$ $R_{4,2}=29\psi_3-50\psi_2,$
$R_{5,2}=84\psi_2^2,$ $R_{3,3}=-12,$ $R_{4,3}=29\psi_2,$  $R_{3,4}=5$ and $R_{l,m}$
is zero otherwise.
Substituting these into~(\ref{e213}), and combining with~(\ref{e212}), we obtain  the following
explicit expression for the left hand side of 
the partial differential equation~(\ref{e21a}):  
\begin{eqnarray} \label{e214}
5760(T_0 G_2-T_1) &=& 
\sum_{i=2}^4 \frac{S_i}{(1-\psi_1)^i} + \sum_{i=3}^5 \frac{2Q_i}{(1-\psi_1)^i}  \nonumber \\
    &\mbox{}& + \sum_{i=4}^7  
     \frac{1}{(1-\psi_1)^i} \sum_m \psi_{m+1} \left( 1+\psi_0-M_{-2,0}-N_0 \right)R_{i-1,m} \nonumber \\ 
&\mbox{}& +  \sum_{i=3}^6   \frac{1}{(1-\psi_1)^i} \sum_m \left(\psi_m-M_{-2,m}-N_m\right)R_{i,m}.
\end{eqnarray}
Then $5760(T_0 G_2-T_1)(1-\psi_1)^7$ is a polynomial of degree 6 in the $q_i$'s, with constant term
equal to zero. Let
$C_i$ be the contribution from terms of (homogeneous) total degree $i$ in the
$q_i$'s for $i=1,\ldots,6.$ The explicit expressions for these are given in Appendix~\ref{A1}.
 
To complete the proof of Theorem~\ref{T1} we show that each of these contributions is identically zero.
It is convenient to introduce the symmetrization operator $\varpi_{1,\ldots,i}$ on the ring
$\cH_i[q_1,q_2,\ldots]$ of homogeneous polynomials of total degree $i$ in $q_1, q_2,\ldots,$
defined by
$$\varpi_{1,\ldots,i} \left(q_{\alpha_1}\cdots q_{\alpha_i} \right)
= \sum_{\pi\in\symgp_i}x_{\pi(1)}^{\alpha_1}\cdots x_{\pi(i)}^{\alpha_i}$$
and extended linearly to the whole of the ring.
Since $\varpi_{1,\ldots,i} f=0$ implies that $f=0$ for $f\in \cH_i[q_1,q_2,\ldots],$ it is sufficient to show that 
$\varpi_{1,\ldots,i} C_i=0$ for $i=1,\ldots,6.$

Now let $w\equiv w(x)$ be the unique solution of the functional equation 
\begin{eqnarray}\label{ew}
w=xe^w,
\end{eqnarray}
and let $w_i\equiv w(x_i),$ where $x_i$ is an indeterminate for $i\ge1.$
Let $w_i^{(j)}= x_i(\partial/\partial x_i)^j w_i.$ Since
$x\partial/\partial x = w/(1-w) \partial/\partial w,$ it is a straightforward
matter to express $w_i^{(k)}$ as a rational function of $w_i.$ 
For example,
\begin{eqnarray*}\label{b2}
\wij1{i} = \frac {w_i}{1-w_i},\quad
\wij2{i} = \frac {w_i}{(1-w_i)^3}, \quad
\wij3{i} = \frac {w_i+2w_i^2}{(1-w_i)^5}.
\end{eqnarray*}
Moreover, 
$w_1, w_2,\ldots$ are algebraically independent.
Now from~\cite{GJdec11}
$$\sum_{i,j,r\ge1} \frac{a_ia_ja_rj^k}{i(i+r)}x_1^{i+j}x_2^r = 
w_1^{(k+2)}\left(
\frac{x_2}{x_1-x_2} - \frac{w_2^{(1)}}{w_1-w_2} -w_2^{(1)}\right),$$
so
\begin{eqnarray*}
\varpi_{1,2}N_k &=& \frac{w_1^{(k+2)}w_2^{(1)}-w_2^{(k+2)}w_1^{(1)}}{w_1-w_2}
+ w_1^{(k+2)}w_2^{(1)}+  w_2^{(k+2)}w_1^{(1)}, \\
\varpi_{1}M_{i,j} &=& w_1^{(i+2)} w_1^{(j+2)}, \\
\varpi_{1} \psi_k &=& w_1^{(k+1)}.
\end{eqnarray*}
Each of these is a polynomial in $W_1$ and $W_2$ where $W_i=1/(1-w_i),$ so 
$\varpi_{1,\ldots,i} C_i$  is a polynomial in $W_1,\ldots,W_i$ alone.
Then $\varpi_{1,\ldots,i}C_i$ may be obtained from the constituents 
$\varpi_{1}\psi_i,$ $\varpi_{1}M_{l,m}$ and $\varpi_{1,2}N_k$  by 
distributing the indeterminates $x_1,\ldots,x_i$ as disjoint subsets of arguments for these
constituents in all possible ways. We have used
{\sf Maple} to carry out this routine but laborious task, and have thus  established that $\varpi_{1,\ldots,i}C_i$
is identically zero as a polynomial in $W_1,\ldots, W_i$  for $i=1,\ldots,6.$ This completes the proof.  \qed

\section{An explicit expression for the number of ramified coverings of the sphere by the double torus}\label{Sexnrc}
In Theorem~\ref{T1} we have determined the generating series $F_2$ for the
ramification numbers $\mu^{(2)}_m(\alpha).$ In this section we expand this series and
thus give an explicit expression for $\mu^{(2)}_m(\alpha).$
The following result is needed, in which $a_\alpha=a_{\alpha_1}a_{\alpha_2}\cdots$
and $m_\nu$ is the monomial symmetric function with exponents specified by the parts of the 
partition $\nu.$

\begin{lemma}\label{L2}
For $\alpha\vdash n,$ $n\ge1,$
$$
\left[x^n p_\alpha\right]\,\frac1{1-\psi_1}\,\prod_{i\ge1}\frac{\psi_i^{j_i}}{j_i!}
=  \frac{a_\alpha}{\macz\alpha} n^{l(\alpha)-l(\nu)}m_\nu(\alpha),
$$
where $\nu=(1^{j_1}\,2^{j_2}\ldots).$
\end{lemma}
\proof  Let
$$\Lambda=\left[x^n \right]\,\frac1{1-\psi_1}\,\prod_{i\ge1}\frac{\psi_i^{j_i}}{j_i!}.$$
Then by Lagrange's Implicit Function Theorem
\begin{eqnarray*}
\Lambda &=& [t^n] \left(\prod_{i\ge1}\frac{\psi_i(t,\bfp)^{j_i}}{j_i!}\right) e^{n\psi_0(t,\bfp)} \\
&=& [t^n]\sum \prod_{i,k\ge1} 
\frac{\left(k^{i-1} a_kp_k t^k\right)^{j_{i,k}}}{j_{i,k}!}\,\,
\prod_{k\ge1}\frac{ \left( nk^{-1}a_kp_kt^k\right)^{d_k}}{d_k!}
\end{eqnarray*}
where the summation is over $j_{i,k}\ge0$, for $i,k\ge1,$ and $d_k\ge0$ for $k\ge1,$ restricted by
$\sum_{k\ge1}j_{i,k}=j_i,$ for $i\ge1.$
Thus, if $\alpha=(1^{b_1}\, 2^{b_2}\,\ldots),$ where $\alpha\vdash n,$ then
$$
[p_\alpha]\Lambda = \frac{a_\alpha}{\macz\alpha} \sum n^{d_1+d_2+\cdots}\prod_{k\ge1}
\frac{b_k!}{d_k!} \prod_{i\ge1} \frac{k^{i j_{i,k}}}{j_{i,k}!},
$$
where the summation is now further restricted by $\sum_{i\ge1}j_{i,k}=b_k-d_k,$ for $k\ge1.$
Then
$$
[p_\alpha]\Lambda = \frac{a_\alpha}{\macz\alpha}
n^{b_1+b_2+\cdots-j_1-j_2-\cdots} \left[y_1^{j_1} y_2^{j_2}\cdots\right]\,
\prod_{k\ge1}\left(1 + \sum_{i\ge1} k^i y_i  \right)^{b_k}
$$
and the result follows. \qed

\medskip

Applying Lemma~\ref{L2} to the generating series $F_2$ given in terms of the $\psi_i$'s in Theorem~\ref{T1},
we immediately obtain the following result which gives an explicit expression for $\mu_m^{(2)}(\alpha).$
This expression is a symmetric function in the parts of $\alpha,$ a linear combination of monomial
symmetric functions. The explicit expression for $\mu_m^{(1)}(\alpha)$ obtained in~\cite{GJdec11}
is a symmetric function expressed in terms of the elementary symmetric functions $e_k(\alpha),$ where $k\ge1.$
These forms are closely related since $e_k=m_{(1^k)}.$


\medskip

\begin{corollary}\label{Cor1}
\begin{eqnarray*}
\mu_m^{(_2)}(\alpha) &=& (n+m+2)! \frac{a_\alpha}{\macz\alpha} \frac{n^m}{5760}
\sum_{k\ge0}\left(
\frac{(k+1)!}{n^{k+1}}
\left( 5 m_{(4\,1^k)} - 12  m_{(3\,1^k) } + 7  m_{(2\,1^k)} \right) \right. \\
& & \mbox{} +  \left.
\frac{(k+2)!}{n^{k+2}}
\left( \frac{29}2 m_{(3\,2\,1^k)} - 25  m_{(2^2\,1^k) }\right) 
+\frac{(k+3)!}{n^{k+3}} 28 m_{(2^3\,1^k)} \right).
\end{eqnarray*}
\end{corollary}
\proof Direct from Theorem~\ref{T1} and Lemma~\ref{L2}. \qed

\section{A proof of the Graber-Pandharipande recurrence equation}\label{SpGPre}
We conclude with an examination of the case $\alpha=(1^n),$ corresponding to
no ramification over $\infty.$ It will be convenient to denote $\mu^{({g})}_n(1^{n})$
by $\mu^{({g})}_n$ for brevity.  For $g\ge0,$ let
$f_g$ be the specialization of $F_g$ with $p_1=1,$ and $p_i=0$ for $i>1.$ Then
$$f_g =\sum_{n\ge1}{x^n}\frac{\mu_n^{(g)}}{(2n+2g-2)!},$$ 
and under these specializations of the $p_i$'s we have $s=w$ where $w$ is the unique solution of the
functional equation~(\ref{ew}) and $\psi_i=w$ for all $i.$ Thus from~(\ref{e4}) and~(\ref{e4a})
we have
\begin{eqnarray*}
\ddx2 f_0 &=& w, \\
f_1 &=& \frac1{24} \left( \log(1-w)^{-1}-w\right),
\end{eqnarray*}
and from Theorem~\ref{T1}
$$f_2=\frac1{5760}\left( \frac{4w^2}{(1-w)^4}+\frac{28w^3}{(1-w)^5}\right).$$
An explicit expression can in fact be obtained for $\mu^{({2})}_n.$
The expression is
\begin{corollary}\label{C0}
$$ \mu^{({2})}_n=\frac{(2n+2)!}{1440 n}\left(12A_4+21A_3+2A_2\right),$$
where
$$A_k=\sum_{i=0}^{n-k}\binom{i+5}{5}\frac{n^{n-i-k}}{(n-i-k)!}.$$
\end{corollary}
\proof The results follows by applying Lagrange's Implicit Function Theorem
to the above expression for $f_2.$ \qed

\medskip
Recurrence equations can be obtained for this number, and our interest in these,
or rather, the corresponding differential equations for $f_2$, is that they
may cast light on a more direct way of obtaining the $\mu^{({2})}_n$.
It is convenient to introduce the operator
$$D=x\frac{d}{dx},$$
and to change variable to $W=1/(1-w).$ Then
$D=W^2(W-1){d}/{dW}.$
Now $Df_0=1/2-W^{-2}/2,$ $D^2f_0=1-1/W,$ and $D^rf_0$ is a polynomial in
$W$ for $r\ge2.$ Moreover, $D^r f_1$ is a polynomial in $W$ for $r\ge 1$
and $D^r f_2$ is a polynomial in $W$ for $r\ge0.$
Then these derivatives are algebraically dependent, so $f_2$ satisfies
a differential equation. Clearly, this equation is not unique.

To decide upon the form that such a differential equation may take we suppose there exists a
(combinatorial or geometric) construction acting on selected sheets that decomposes a
covering into two connected coverings whose genera sum to the genus of the original covering.
The  combinatorial effect of $D$ is to select a single sheet in all possible ways.
We therefore seek a formal linear differential equation for $f_2$
that involves terms of the form $(D^p f_i)(D^q f_j)$ where $i+j=2,$
with the additional condition that it suffices to select at most four sheets, so $p+q\le4$
and $p,q\ge1,$ together with terms of the form $D^rf_1$ where $2\le r\le 3.$
Such a differential equation has the form
\begin{eqnarray*}
\left( b_1 D^2+b_2D+b_3\right)f_2 &=& \left(b_4 D^3+b_5D^2\right)f_1
+b_6\left(D^2f_0\right) \left(D^2 f_2\right)
+b_7 \left(D^2 f_1\right)^2
+b_8 \left(Df_1\right)\left(D^3 f_1\right) \\
&\mbox{}&
+b_9\left(D^2f_2\right)\left(Df_0\right)
+b_{10}\left(D^2f_0\right)\left(Df_2\right).
\end{eqnarray*}
It follows by substituting the computed derivatives into the differential
equation, equating coefficients of powers of $W$ to obtain a system of homogeneous
linear equations and solving this system, that the solution space is 4 dimensional
and that
\begin{eqnarray*}
b_3 &=& -4 b_1 - 2 b_2 + 240 b_4 + 120 b_5, \\
b_6 &=& - \frac{11}{2} b_1 - \frac{3}{2} b_2 - 72 b_4 - 70 b_5, \\
b_7 &=& \frac{47}4 b_1 + \frac{23}4 b_2 - 1236 b_4 - 875 b_5, \\
b_8 &=& - \frac{293}4 b_1 - \frac{85}4 b_2 - 264 b_4 - 420 b_5, \\
b_9 &=& 13 b_1 + 3 b_2 + 144 b_4 + 140 b_5, \\
b_{10} &=& \frac{35}2 b_1 + \frac{7}2 b_2 + 336 b_4 + 280 b_5, \\
\end{eqnarray*}
where $b_1, b_2, b_4, b_5$ are arbitrary. The system therefore has nontrivial solutions.
\begin{corollary}\label{C1}
\begin{eqnarray*}
\mu^{({2})}_n &=&  n^2\left(\frac{97}{136}n-\frac{20}{17}  \right) \mu^{({1})}_n
+ \sum_{j=1}^{n-1}  \binom{2n}{2j-2} 
  \mu^{({0})}_j\,\mu^{({2})}_{n-j} j(n-j)\left(8n-\frac{115}{17}j \right)\nonumber\\
& & \mbox{} + \sum_{j=1}^{n-1}  \binom{2n}{2j} \mu^{({1})}_j\,
\mu^{({1})}_{n-j} j(n-j)
\left(\frac{11697}{34}j(n-j)-\frac{3899}{68}n^2 \right).
\end{eqnarray*}
\end{corollary}
\proof
By setting $ b_1=4,$ $ b_2= 6,$ $ b_4= 97/136,$ $b_5= -20/17$
we have
\begin{eqnarray*}
\left(4\ddx 2 +6\dx +2\right)f_2 &=& \left( \frac{97}{136}\ddx3-\frac{20}{17}\ddx2 \right)f_1
+ 8 \left( \ddx2 f_2 \right) \left( \dx f_0 \right)
+ \frac{21}{17} \left( \ddx2 f_0 \right) \left( \dx f_2 \right)
+ \frac{3899}{17}\left( \ddx2 f_1\right)^2 \nonumber \\
& & \mbox{} -\frac{3899}{34}\left(\dx f_1\right) \left( \ddx3 f_1\right),
\end{eqnarray*}
and the result follows immediately. \qed

\medskip
This establishes the recurrence equation for $\mu_n^{(2)}(1^n),$
corresponding to the case of simple ramification,
conjectured by Graber and Pandharipande~\cite{V1}.

\section{A ``simpler'' relationship for the double torus}\label{samrdt}
In view of the combinatorial interpretation of the differential operator $D,$ a
differential equation such as the one given in the proof of Corollary~\ref{C1}
(or, equivalently, a recurrence equation) may have a more direct
combinatorial explanation, which may in turn suggest a geometrical explanation.
For this purpose it is therefore prudent to look for a differential equation with fewer
terms, and  whose coefficients are, at the very least, potentially more susceptible to
combinatorial explanation. Since there are four independent parameters $b_1, b_2, b_4, b_5$,
we may impose three further conditions to lessen the number of terms in the
differential equation, and then divide out the remaining parameter.
The obvious conditions to apply are those that remove terms from the differential equation.
The next two corollaries give instances where these criteria are met.

The first instance is a second order linear differential equation for $f_2$ with 
simple coefficients that has contributions from the sphere, torus and the
double torus on the right hand side.

\begin{corollary}\label{C2}
\begin{eqnarray*}
\left(2\ddx 2 -6\dx +2\right)f_2 = \left( \frac{1}{24}\ddx3-\frac{1}{10}\ddx2 \right)f_1
+ 2 \left( \ddx2 f_0 \right) \left( \ddx2 f_2 \right)
+ 25\left( \ddx2 f_1\right)^2
+12\left(\dx f_1\right) \left( \ddx3 f_1\right).
\end{eqnarray*}
\end{corollary}
\proof Set $b_9=b_{10}=0,$ and $b_1=2.$ \qed

\medskip
The second instance is obtained by imposing conditions
to eliminate the presence of contributions from the double torus (and therefore
the sphere) on the right hand side of the differential equation.
This gives a first order linear differential equation for $f_2$ with  simple
coefficients that has contributions only from the torus on the right hand side.

\begin{corollary}\label{C3}
\begin{eqnarray*}
\left(2\dx +3\right)f_2 = \frac14\left( \frac{5}{12}\ddx3-\frac{3}{5}\ddx2 \right)f_1
+ 14 \left( \ddx2 f_1\right)^2
-7\left(\dx f_1\right) \left( \ddx3 f_1\right).
\end{eqnarray*}
\end{corollary}
\proof Set $b_6=b_9=b_{10}=0.$ \qed



\medskip
It remains to determine whether a differential equation can be found
with fewer terms and with equally simple coefficients.
To assist in the search, let $\Delta f$, for a polynomial $f$ in $W$,
be defined to be $(k,l)$, where $k$ and $l$ are, respectively, the lowest and highest degrees
of the terms of $f$ in $W.$ We will refer to this as the {\em degree span} of $f.$
Then it is readily seen that
$\Delta D^rf_0=(r-2,2r-5)$ for $r\ge4,$ 
$\Delta D^3 f_0=(0,1),$
$\Delta D^2 f_0=(-1,0)$ and
$\Delta D f_0=(-2,0).$
Also $\Delta D f_1=(0,2)$,
$\Delta D^rf_1=(r,2r)$ for $r\ge2,$ and $\Delta D^rf_2=(r+2,2r+5)$ for $r\ge0.$ 

Consider $(b_1 D+ b_2)f_2$ where $b_1$ and $b_2$ are generic real numbers.
Then, from the above spans, $\Delta(b_1 D+ b_2)f_2=(2,7).$
We now construct another expression from $f_1$ and $f_0$ that has the same span.
From the above expressions for spans,
$\Delta(D^2f_1)(Df_1)=(2,6),$
$\Delta D^2f_1=(2,4),$ 
$\Delta D^3f_1=(3,6).$
and $\Delta D^6 f_0=(4,7).$
Then, for generic $b_3,\ldots,b_6,$ we have
$\Delta(b_3D^6f_0 + b_4(D^2f_1)(Df_1) + b_5D^2f_1 + b_6D^3f_1)=(2,7).$
Therefore we consider the differential equation
$$(b_1 D+ b_2)f_2 = b_3D^6f_0 + b_4(D^2f_1)(Df_1) + b_5D^2f_1 + b_6D^3f_1.$$
\begin{corollary}\label{C4}
$$f_2= \frac{1}{6!} (7D^3 - 8D^2)f_1 -\frac{14}{15}(D^2f_1)(Df_1).$$
Moreover,
$$ \mu_n^{(2)} =
2\binom{2n+2}{2}
\left( \frac{1}{6!}(7n^3-8n^2)\mu_n^{(1)}
-\frac{7n}{15}\sum_{j=1}^{n-1}j(n-j)\binom{2n}{2j}\mu_j^{(1)}\mu_{n-j}^{(1)}\right).$$
\end{corollary}
\proof
By the argument of the previous section we obtain a solution space of dimension $1,$
for the system of linear equations, and this gives a unique differential equation
up to a normalizing factor. The recurrence equation is obtained by comparing
coefficients of $x^n$ on each side of the equation. \qed

\medskip
The above corollary gives an equation for $f_2$ that certainly has fewer terms
than the differential equations of the earlier corollaries.
Morover, $7=2^3-1,$ $8=2^3$ and $15=3!!,$
(where $n!!=(2n)!/2^nn!$, the number of perfect matchings on a $2n$-set) each of which is a 
number with a known combinatorial interpretation.


Implicit in the above discussion is the assumption that there
is no linear recurrence equation for $\mu_n^{(2)}$. This can be verified
easily through an algebraic argument that appeals to the fact that $e^w$ is a
transcendental series in $w.$


\section*{Acknowledgements}
This work was supported by grants individually
to IPG and DMJ from the Natural
Sciences and Engineering Research Council of Canada.



\appendix

\section{Ramified coverings of the sphere by the triple torus}\label{A0}
\begin{eqnarray*}
F_3(x,\bfp) &=&\frac{1}{2^39!}  
\,{\displaystyle \frac {35\,{\psi_{7}} - 147\,{\psi_{6}} + 205\,{\psi_{5}} 
- 93\,{\psi_{4}}}{(1 - {\psi_{1}})^{5}}} \\ 
& & + {\displaystyle \frac {1}{2^39!}} (- 930\,{\psi_{2}}\,{\psi_{3}} + 607\,{\psi_{4}}^{2} 
+ 1501\,{\psi_{3}}^{2} + 2329\,{\psi_{2}}\,{\psi_{4}} + 539\,{\psi_{2}}\,{\psi_{6}} 
 + 1006\,{\psi_{3}}\,{\psi_{5}} \\
& & \mbox{} - 3078\,{\psi_{3}}\,{\psi_{4}} - 1938\,{\psi_{2}}\,{\psi_{5}})
 \left/ {\vrule height0.54em width0em depth0.54em} \right. \! 
 \! (1 - {\psi_{1}})^{6}\mbox{} \\
& &+ {\displaystyle \frac {1}{2^39!}} ( 
 13452\,{\psi_{2}}\,{\psi_{3}}\,{\psi_{4}} + 2915
\,{\psi_{3}}^{3} - 16821\,{\psi_{2}}\,{\psi_{3}}^{2}
 - 12984\,{\psi_{2}}^{2}\,{\psi_{4}} \\
 & & \mbox{} + 12885\,{\psi_{2}}^{2}\,{\psi_{3}} + 4284\,
{\psi_{2}}^{2}\,{\psi_{5}} - 1395\,{\psi_{2}}^{3})
 \left/ {\vrule height0.54em width0em depth0.54em} \right. \! 
 \! (1 - {\psi_{1}})^{7} \\
 & & \mbox{} + {\displaystyle \frac {1}{2^39!}} \,
{\displaystyle \frac {22260\,{\psi_{2}}^{3}\,{\psi_{4}}
 + 43050\,{\psi_{2}}^{2}\,{\psi_{3}}^{2} - 55300\,{\psi_{2}}^{3}\,{\psi_{3}} + 10710\,{\psi_{2}}^{4}}{(1 - {
\psi_{1}})^{8}}}  \\
 & & \mbox{} + {\displaystyle \frac {1}{2^39!}} \,
{\displaystyle \frac {81060\,{\psi_{2}}^{4}\,{\psi_{3}}
 - 31220\,{\psi_{2}}^{5}}{(1 - {\psi_{1}})^{9}}}  + 
{\displaystyle \frac {245}{20736}} \,{\displaystyle \frac {{\psi_{2}}^{6}}{(1 - {\psi_{1}})^{10}}} 
\end{eqnarray*}


\section{The expression for $5760(T_0G_2-T_1)(1-\psi_1)^7$}\label{A1}
Listed below for $i=1,\ldots,6$ is $C_i,$  the contribution to $5760(T_0G_2-T_1)(1-\psi_1)^7$ from terms
of total degree $i$ in the $q_j$'s. The $M_{l,m},$ $N_k$ and $\psi_i$
are series in the $q_j$'s. 

\begin{eqnarray*}
C_1 &=& \lefteqn{ - 240\,{M_{2, \,0}} - 125\,{M_{1, \,1}} + 250\,{M_{1, 
\,0}} - 5\,{M_{0, \,0}} - 7\,{M_{-2, \,2}} + 3\,{\psi _{4}} - 29
\,{\psi _{3}} + 21\,{\psi _{2}} + 12\,{M_{-2, \,3}}} \\
 & & \mbox{} - 5\,{M_{-2, \,4}} + 5\,{\psi _{5}}
\mbox{\hspace{336pt}}
\end{eqnarray*}

\begin{eqnarray*}
C_2 &=&\lefteqn{20\,{\psi _{1}}\,{M_{-2, \,4}} + 29\,{\psi _{3}}^{2} - 9
\,{\psi _{1}}\,{\psi _{4}} + 87\,{\psi _{1}}\,{\psi _{3}} - 63\,{
\psi _{1}}\,{\psi _{2}} - 48\,{\psi _{1}}\,{M_{-2, \,3}} - 490\,{
\psi _{2}}\,{M_{1, \,0}}} \\
 & & \mbox{} + 130\,{\psi _{2}}\,{M_{0, \,0}} + 28\,{\psi _{1}}\,
{M_{-2, \,2}} + 7\,{\psi _{3}}\,{\psi _{0}} - 7\,{\psi _{3}}\,{M
_{-2, \,0}} - 29\,{\psi _{2}}\,{M_{-2, \,3}} + 5\,{\psi _{5}}\,{
\psi _{0}} - 5\,{\psi _{5}}\,{M_{-2, \,0}} \\
 & & \mbox{} - 12\,{\psi _{4}}\,{\psi _{0}} + 12\,{\psi _{4}}\,{M
_{-2, \,0}} - 15\,{\psi _{1}}\,{\psi _{5}} - 29\,{M_{-2, \,2}}\,{
\psi _{3}} + 50\,{M_{-2, \,2}}\,{\psi _{2}} - 15\,{M_{-2, \,1}}\,
{\psi _{4}} \\
 & & \mbox{} + 36\,{M_{-2, \,1}}\,{\psi _{3}} - 21\,{M_{-2, \,1}}
\,{\psi _{2}} + 625\,{\psi _{1}}\,{M_{1, \,1}} - 1250\,{\psi _{1}
}\,{M_{1, \,0}} + 25\,{\psi _{1}}\,{M_{0, \,0}} - 5\,{N_{4}} + 12
\,{N_{3}} \\
 & & \mbox{} - 7\,{N_{2}} - 120\,{\psi _{3}}\,{M_{0, \,0}} + 30\,
{\psi _{2}}\,{\psi _{3}} - 79\,{\psi _{2}}^{2} + 44\,{\psi _{4}}
\,{\psi _{2}} + 1200\,{\psi _{1}}\,{M_{2, \,0}}
\end{eqnarray*}

\begin{eqnarray*}
C_3 &=&\lefteqn{15\,{\psi _{1}}^{2}\,{\psi _{5}} + 36\,{\psi _{1}}\,{
\psi _{4}}\,{\psi _{0}} - 36\,{\psi _{1}}\,{\psi _{4}}\,{M_{-2, 
\,0}} - 15\,{\psi _{1}}\,{\psi _{5}}\,{\psi _{0}} + 15\,{\psi _{1
}}\,{\psi _{5}}\,{M_{-2, \,0}}} \\
 & & \mbox{} + 45\,{\psi _{1}}\,{M_{-2, \,1}}\,{\psi _{4}} - 108
\,{\psi _{1}}\,{M_{-2, \,1}}\,{\psi _{3}} + 63\,{\psi _{1}}\,{M_{
-2, \,1}}\,{\psi _{2}} + 87\,{\psi _{1}}\,{M_{-2, \,2}}\,{\psi _{
3}} + 20\,{\psi _{1}}\,{N_{4}} \\
 & & \mbox{} - 42\,{\psi _{1}}^{2}\,{M_{-2, \,2}} + 9\,{\psi _{1}
}^{2}\,{\psi _{4}} - 87\,{\psi _{1}}^{2}\,{\psi _{3}} + 63\,{\psi
 _{1}}^{2}\,{\psi _{2}} + 72\,{\psi _{1}}^{2}\,{M_{-2, \,3}} - 30
\,{\psi _{1}}^{2}\,{M_{-2, \,4}} \\
 & & \mbox{} - 29\,{N_{2}}\,{\psi _{3}} + 50\,{N_{2}}\,{\psi _{2}
} - 15\,{N_{1}}\,{\psi _{4}} + 36\,{N_{1}}\,{\psi _{3}} - 21\,{N
_{1}}\,{\psi _{2}} + 1960\,{\psi _{1}}\,{\psi _{2}}\,{M_{1, \,0}}
 - 48\,{\psi _{1}}\,{N_{3}} \\
 & & \mbox{} + 28\,{\psi _{1}}\,{N_{2}} - 520\,{\psi _{1}}\,{\psi
 _{2}}\,{M_{0, \,0}} + 480\,{\psi _{1}}\,{\psi _{3}}\,{M_{0, \,0}
} - 88\,{\psi _{1}}\,{\psi _{4}}\,{\psi _{2}} - 29\,{\psi _{3}}^{
2}\,{M_{-2, \,0}} - 58\,{\psi _{1}}\,{\psi _{3}}^{2} \\
 & & \mbox{} - 116\,{M_{-2, \,1}}\,{\psi _{2}}\,{\psi _{3}} + 86
\,{\psi _{2}}\,{\psi _{3}}\,{M_{-2, \,0}} - 60\,{\psi _{1}}\,{
\psi _{2}}\,{\psi _{3}} + 44\,{\psi _{2}}\,{\psi _{4}}\,{\psi _{0
}} - 86\,{\psi _{2}}\,{\psi _{3}}\,{\psi _{0}} \\
 & & \mbox{} - 44\,{\psi _{2}}\,{\psi _{4}}\,{M_{-2, \,0}} + 21\,
{\psi _{0}}\,{\psi _{2}}^{2} - 21\,{M_{-2, \,0}}\,{\psi _{2}}^{2}
 + 158\,{\psi _{1}}\,{\psi _{2}}^{2} + 100\,{M_{-2, \,1}}\,{\psi 
_{2}}^{2} - 84\,{\psi _{2}}^{2}\,{M_{-2, \,2}} \\
 & & \mbox{} + 29\,{\psi _{3}}^{2}\,{\psi _{0}} - 1250\,{\psi _{1
}}^{2}\,{M_{1, \,1}} + 2500\,{\psi _{1}}^{2}\,{M_{1, \,0}} - 50\,
{\psi _{1}}^{2}\,{M_{0, \,0}} - 2400\,{\psi _{1}}^{2}\,{M_{2, \,0
}} \\
 & & \mbox{} - 21\,{\psi _{1}}\,{\psi _{3}}\,{\psi _{0}} - 150\,{
\psi _{1}}\,{M_{-2, \,2}}\,{\psi _{2}} + 21\,{\psi _{1}}\,{\psi 
_{3}}\,{M_{-2, \,0}} + 87\,{\psi _{1}}\,{\psi _{2}}\,{M_{-2, \,3}
} + 40\,{\psi _{2}}^{3} - 29\,{N_{3}}\,{\psi _{2}} \\
 & & \mbox{} - 245\,{\psi _{2}}^{2}\,{M_{0, \,0}} - 5\,{\psi _{5}
}\,{N_{0}} - 7\,{\psi _{3}}\,{N_{0}} + 12\,{\psi _{4}}\,{N_{0}}
 + 200\,{\psi _{3}}\,{\psi _{2}}^{2}
\end{eqnarray*}

\begin{eqnarray*}
C_4 &=&\lefteqn{200\,{\psi _{2}}^{2}\,{\psi _{3}}\,{\psi _{0}} - 200\,{
\psi _{2}}^{2}\,{\psi _{3}}\,{M_{-2, \,0}} - 100\,{\psi _{0}}\,{
\psi _{2}}^{3} - 3\,{\psi _{1}}^{3}\,{\psi _{4}} + 29\,{\psi _{1}
}^{3}\,{\psi _{3}} - 21\,{\psi _{1}}^{3}\,{\psi _{2}}} \\
 & & \mbox{} - 48\,{\psi _{1}}^{3}\,{M_{-2, \,3}} + 20\,{\psi _{1
}}^{3}\,{M_{-2, \,4}} - 140\,{\psi _{2}}^{3}\,{M_{-2, \,1}} + 100
\,{M_{-2, \,0}}\,{\psi _{2}}^{3} + 44\,{\psi _{1}}^{2}\,{\psi _{4
}}\,{\psi _{2}} \\
 & & \mbox{} - 2940\,{\psi _{1}}^{2}\,{\psi _{2}}\,{M_{1, \,0}}
 + 780\,{\psi _{1}}^{2}\,{\psi _{2}}\,{M_{0, \,0}} - 720\,{\psi 
_{1}}^{2}\,{\psi _{3}}\,{M_{0, \,0}} + 28\,{\psi _{1}}^{3}\,{M_{
-2, \,2}} \\
 & & \mbox{} + 172\,{\psi _{1}}\,{\psi _{2}}\,{\psi _{3}}\,{\psi 
_{0}} + 88\,{\psi _{1}}\,{\psi _{2}}\,{\psi _{4}}\,{M_{-2, \,0}}
 - 172\,{\psi _{1}}\,{\psi _{2}}\,{\psi _{3}}\,{M_{-2, \,0}} + 
232\,{\psi _{1}}\,{M_{-2, \,1}}\,{\psi _{2}}\,{\psi _{3}} \\
 & & \mbox{} - 200\,{\psi _{1}}\,{M_{-2, \,1}}\,{\psi _{2}}^{2}
 + 29\,{\psi _{1}}^{2}\,{\psi _{3}}^{2} + 72\,{\psi _{1}}^{2}\,{N
_{3}} - 42\,{\psi _{1}}^{2}\,{N_{2}} - 30\,{\psi _{1}}^{2}\,{N_{4
}} - 88\,{\psi _{1}}\,{\psi _{2}}\,{\psi _{4}}\,{\psi _{0}} \\
 & & \mbox{} - 58\,{\psi _{1}}\,{\psi _{3}}^{2}\,{\psi _{0}} + 58
\,{\psi _{1}}\,{\psi _{3}}^{2}\,{M_{-2, \,0}} + 30\,{\psi _{1}}^{
2}\,{\psi _{2}}\,{\psi _{3}} - 42\,{\psi _{1}}\,{\psi _{0}}\,{
\psi _{2}}^{2} + 42\,{\psi _{1}}\,{M_{-2, \,0}}\,{\psi _{2}}^{2}
 \\
 & & \mbox{} + 86\,{\psi _{3}}\,{N_{0}}\,{\psi _{2}} - 79\,{\psi 
_{1}}^{2}\,{\psi _{2}}^{2} + 100\,{N_{1}}\,{\psi _{2}}^{2} - 29\,
{\psi _{3}}^{2}\,{N_{0}} - 21\,{\psi _{2}}^{2}\,{N_{0}} - 40\,{
\psi _{1}}\,{\psi _{2}}^{3} \\
 & & \mbox{} - 116\,{N_{1}}\,{\psi _{2}}\,{\psi _{3}} - 108\,{
\psi _{1}}\,{N_{1}}\,{\psi _{3}} + 63\,{\psi _{1}}\,{N_{1}}\,{
\psi _{2}} + 735\,{\psi _{1}}\,{\psi _{2}}^{2}\,{M_{0, \,0}} - 36
\,{\psi _{1}}\,{\psi _{4}}\,{N_{0}} \\
 & & \mbox{} + 21\,{\psi _{1}}\,{\psi _{3}}\,{N_{0}} + 87\,{\psi 
_{1}}\,{N_{3}}\,{\psi _{2}} + 15\,{\psi _{1}}\,{\psi _{5}}\,{N_{0
}} - 45\,{\psi _{1}}^{2}\,{M_{-2, \,1}}\,{\psi _{4}} + 87\,{\psi 
_{1}}\,{N_{2}}\,{\psi _{3}} \\
 & & \mbox{} - 150\,{\psi _{1}}\,{N_{2}}\,{\psi _{2}} + 45\,{\psi
 _{1}}\,{N_{1}}\,{\psi _{4}} - 200\,{\psi _{1}}\,{\psi _{3}}\,{
\psi _{2}}^{2} - 21\,{\psi _{1}}^{2}\,{\psi _{3}}\,{M_{-2, \,0}}
 - 87\,{\psi _{1}}^{2}\,{\psi _{2}}\,{M_{-2, \,3}} \\
 & & \mbox{} + 15\,{\psi _{1}}^{2}\,{\psi _{5}}\,{\psi _{0}} - 15
\,{\psi _{1}}^{2}\,{\psi _{5}}\,{M_{-2, \,0}} - 36\,{\psi _{1}}^{
2}\,{\psi _{4}}\,{\psi _{0}} + 36\,{\psi _{1}}^{2}\,{\psi _{4}}\,
{M_{-2, \,0}} + 108\,{\psi _{1}}^{2}\,{M_{-2, \,1}}\,{\psi _{3}}
 \\
 & & \mbox{} - 87\,{\psi _{1}}^{2}\,{M_{-2, \,2}}\,{\psi _{3}} + 
150\,{\psi _{1}}^{2}\,{M_{-2, \,2}}\,{\psi _{2}} + 21\,{\psi _{1}
}^{2}\,{\psi _{3}}\,{\psi _{0}} - 63\,{\psi _{1}}^{2}\,{M_{-2, \,
1}}\,{\psi _{2}} - 5\,{\psi _{1}}^{3}\,{\psi _{5}} \\
 & & \mbox{} + 168\,{\psi _{1}}\,{\psi _{2}}^{2}\,{M_{-2, \,2}}
 + 2400\,{\psi _{1}}^{3}\,{M_{2, \,0}} + 1250\,{\psi _{1}}^{3}\,{
M_{1, \,1}} - 2500\,{\psi _{1}}^{3}\,{M_{1, \,0}} + 50\,{\psi _{1
}}^{3}\,{M_{0, \,0}} \\
 & & \mbox{} - 84\,{N_{2}}\,{\psi _{2}}^{2} - 44\,{\psi _{4}}\,{N
_{0}}\,{\psi _{2}} + 140\,{\psi _{2}}^{4}
\end{eqnarray*}

\begin{eqnarray*}
C_5 &=& \lefteqn{100\,{\psi _{2}}^{3}\,{N_{0}} + 200\,{\psi _{1}}\,{\psi 
_{2}}^{2}\,{\psi _{3}}\,{M_{-2, \,0}} + 1960\,{\psi _{1}}^{3}\,{
\psi _{2}}\,{M_{1, \,0}} + 140\,{\psi _{1}}\,{\psi _{2}}^{3}\,{M
_{-2, \,1}}} \\
 & & \mbox{} - 200\,{\psi _{1}}\,{\psi _{2}}^{2}\,{\psi _{3}}\,{
\psi _{0}} + 140\,{\psi _{2}}^{4}\,{\psi _{0}} - 140\,{\psi _{2}}
^{4}\,{M_{-2, \,0}} - 48\,{\psi _{1}}^{3}\,{N_{3}} + 28\,{\psi _{
1}}^{3}\,{N_{2}} + 20\,{\psi _{1}}^{3}\,{N_{4}} \\
 & & \mbox{} - 7\,{\psi _{1}}^{4}\,{M_{-2, \,2}} + 12\,{\psi _{1}
}^{4}\,{M_{-2, \,3}} - 5\,{\psi _{1}}^{4}\,{M_{-2, \,4}} - 12\,{
\psi _{1}}^{3}\,{\psi _{4}}\,{M_{-2, \,0}} + 100\,{\psi _{1}}\,{
\psi _{0}}\,{\psi _{2}}^{3} \\
 & & \mbox{} - 100\,{\psi _{1}}\,{M_{-2, \,0}}\,{\psi _{2}}^{3}
 - 520\,{\psi _{1}}^{3}\,{\psi _{2}}\,{M_{0, \,0}} + 480\,{\psi 
_{1}}^{3}\,{\psi _{3}}\,{M_{0, \,0}} - 36\,{\psi _{1}}^{3}\,{M_{
-2, \,1}}\,{\psi _{3}} \\
 & & \mbox{} + 21\,{\psi _{1}}^{3}\,{M_{-2, \,1}}\,{\psi _{2}} + 
29\,{\psi _{1}}^{3}\,{M_{-2, \,2}}\,{\psi _{3}} - 50\,{\psi _{1}}
^{3}\,{M_{-2, \,2}}\,{\psi _{2}} - 7\,{\psi _{1}}^{3}\,{\psi _{3}
}\,{\psi _{0}} + 7\,{\psi _{1}}^{3}\,{\psi _{3}}\,{M_{-2, \,0}}
 \\
 & & \mbox{} + 29\,{\psi _{1}}^{3}\,{\psi _{2}}\,{M_{-2, \,3}} - 
5\,{\psi _{1}}^{3}\,{\psi _{5}}\,{\psi _{0}} + 5\,{\psi _{1}}^{3}
\,{\psi _{5}}\,{M_{-2, \,0}} + 12\,{\psi _{1}}^{3}\,{\psi _{4}}\,
{\psi _{0}} - 735\,{\psi _{1}}^{2}\,{\psi _{2}}^{2}\,{M_{0, \,0}}
 \\
 & & \mbox{} - 21\,{\psi _{1}}^{2}\,{\psi _{3}}\,{N_{0}} - 87\,{
\psi _{1}}^{2}\,{N_{3}}\,{\psi _{2}} - 15\,{\psi _{1}}^{2}\,{\psi
 _{5}}\,{N_{0}} + 15\,{\psi _{1}}^{3}\,{M_{-2, \,1}}\,{\psi _{4}}
 + 36\,{\psi _{1}}^{2}\,{\psi _{4}}\,{N_{0}} \\
 & & \mbox{} - 21\,{\psi _{1}}^{2}\,{M_{-2, \,0}}\,{\psi _{2}}^{2
} + 100\,{\psi _{1}}^{2}\,{M_{-2, \,1}}\,{\psi _{2}}^{2} - 84\,{
\psi _{1}}^{2}\,{\psi _{2}}^{2}\,{M_{-2, \,2}} + 29\,{\psi _{1}}
^{2}\,{\psi _{3}}^{2}\,{\psi _{0}} \\
 & & \mbox{} - 29\,{\psi _{1}}^{2}\,{\psi _{3}}^{2}\,{M_{-2, \,0}
} + 86\,{\psi _{1}}^{2}\,{\psi _{2}}\,{\psi _{3}}\,{M_{-2, \,0}}
 - 116\,{\psi _{1}}^{2}\,{M_{-2, \,1}}\,{\psi _{2}}\,{\psi _{3}}
 - 87\,{\psi _{1}}^{2}\,{N_{2}}\,{\psi _{3}} \\
 & & \mbox{} + 150\,{\psi _{1}}^{2}\,{N_{2}}\,{\psi _{2}} - 45\,{
\psi _{1}}^{2}\,{N_{1}}\,{\psi _{4}} + 108\,{\psi _{1}}^{2}\,{N_{
1}}\,{\psi _{3}} - 63\,{\psi _{1}}^{2}\,{N_{1}}\,{\psi _{2}} - 86
\,{\psi _{1}}^{2}\,{\psi _{2}}\,{\psi _{3}}\,{\psi _{0}} \\
 & & \mbox{} - 44\,{\psi _{1}}^{2}\,{\psi _{2}}\,{\psi _{4}}\,{M
_{-2, \,0}} - 172\,{\psi _{1}}\,{\psi _{3}}\,{N_{0}}\,{\psi _{2}}
 + 88\,{\psi _{1}}\,{\psi _{4}}\,{N_{0}}\,{\psi _{2}} + 44\,{\psi
 _{1}}^{2}\,{\psi _{2}}\,{\psi _{4}}\,{\psi _{0}} \\
 & & \mbox{} + 232\,{\psi _{1}}\,{N_{1}}\,{\psi _{2}}\,{\psi _{3}
} - 200\,{\psi _{1}}\,{N_{1}}\,{\psi _{2}}^{2} + 58\,{\psi _{1}}
\,{\psi _{3}}^{2}\,{N_{0}} + 42\,{\psi _{1}}\,{\psi _{2}}^{2}\,{N
_{0}} + 168\,{\psi _{1}}\,{N_{2}}\,{\psi _{2}}^{2} \\
 & & \mbox{} + 21\,{\psi _{1}}^{2}\,{\psi _{0}}\,{\psi _{2}}^{2}
 - 1200\,{\psi _{1}}^{4}\,{M_{2, \,0}} - 625\,{\psi _{1}}^{4}\,{M
_{1, \,1}} + 1250\,{\psi _{1}}^{4}\,{M_{1, \,0}} - 25\,{\psi _{1}
}^{4}\,{M_{0, \,0}} \\
 & & \mbox{} - 140\,{N_{1}}\,{\psi _{2}}^{3} - 200\,{\psi _{3}}\,
{N_{0}}\,{\psi _{2}}^{2}
\end{eqnarray*}

\begin{eqnarray*}
C_6 &=& \lefteqn{5\,{\psi _{1}}^{3}\,{\psi _{5}}\,{N_{0}} - 12\,{\psi _{1
}}^{3}\,{\psi _{4}}\,{N_{0}} + 7\,{\psi _{1}}^{3}\,{\psi _{3}}\,{
N_{0}} + 29\,{\psi _{1}}^{3}\,{N_{3}}\,{\psi _{2}} - 116\,{\psi 
_{1}}^{2}\,{N_{1}}\,{\psi _{2}}\,{\psi _{3}}} \\
 & & \mbox{} + 86\,{\psi _{1}}^{2}\,{\psi _{3}}\,{N_{0}}\,{\psi 
_{2}} - 44\,{\psi _{1}}^{2}\,{\psi _{4}}\,{N_{0}}\,{\psi _{2}} + 
240\,{\psi _{1}}^{5}\,{M_{2, \,0}} + 125\,{\psi _{1}}^{5}\,{M_{1
, \,1}} - 250\,{\psi _{1}}^{5}\,{M_{1, \,0}} \\
 & & \mbox{} + 5\,{\psi _{1}}^{5}\,{M_{0, \,0}} + 100\,{\psi _{1}
}^{2}\,{N_{1}}\,{\psi _{2}}^{2} - 29\,{\psi _{1}}^{2}\,{\psi _{3}
}^{2}\,{N_{0}} - 21\,{\psi _{1}}^{2}\,{\psi _{2}}^{2}\,{N_{0}} - 
84\,{\psi _{1}}^{2}\,{N_{2}}\,{\psi _{2}}^{2} \\
 & & \mbox{} + 29\,{\psi _{1}}^{3}\,{N_{2}}\,{\psi _{3}} - 50\,{
\psi _{1}}^{3}\,{N_{2}}\,{\psi _{2}} + 15\,{\psi _{1}}^{3}\,{N_{1
}}\,{\psi _{4}} - 36\,{\psi _{1}}^{3}\,{N_{1}}\,{\psi _{3}} + 21
\,{\psi _{1}}^{3}\,{N_{1}}\,{\psi _{2}} \\
 & & \mbox{} + 245\,{\psi _{1}}^{3}\,{\psi _{2}}^{2}\,{M_{0, \,0}
} + 12\,{\psi _{1}}^{4}\,{N_{3}} - 7\,{\psi _{1}}^{4}\,{N_{2}} - 
5\,{\psi _{1}}^{4}\,{N_{4}} + 200\,{\psi _{1}}\,{\psi _{3}}\,{N_{
0}}\,{\psi _{2}}^{2} \\
 & & \mbox{} - 100\,{\psi _{1}}\,{\psi _{2}}^{3}\,{N_{0}} + 140\,
{\psi _{1}}\,{N_{1}}\,{\psi _{2}}^{3} - 490\,{\psi _{1}}^{4}\,{
\psi _{2}}\,{M_{1, \,0}} + 130\,{\psi _{1}}^{4}\,{\psi _{2}}\,{M
_{0, \,0}} - 120\,{\psi _{1}}^{4}\,{\psi _{3}}\,{M_{0, \,0}} \\
 & & \mbox{} - 140\,{\psi _{2}}^{4}\,{N_{0}}
\end{eqnarray*}

\end{document}